\documentclass[s5,11pt]{article}
\newtheorem{theorem}{Theorem} 
\newtheorem {corollary}{Corollary}
\newtheorem{lemma}{Lemma}
 
\author{Tord Sj\"odin, Ume\aa} 
\title{ A Note on Gram--Schmidt's   Algorithm for a General Angle.}
\begin{document}%generates the title
    
\maketitle 
 
 \begin{abstract}  The Gram--Schmidt algorithm produces a pairwise orthogonal set $\{ y_i\}^n_1$ from a linearly independent set $\{ x_i\}^n_1$ of vectors in an inner product vector space $V$. We give a linear algorithm that constructs vectors $\{ z_i\}^n_1$ with the same span as $\{ x_i\}^n_1$ and which have pairwise the same prescribed angle or distance, in all cases where this is possible. Finally, we prove an asymptotic property of $\{ z_i\}^\infty_1$ in the infinite dimensional case.\end{abstract}

 
\paragraph{ \it AMS 2000 Subject Classsification:} Primary 15 A 03, Secondary 
15 A 63, 46 C 05

\section { Introduction.}  Let $V$ be a real vector space with inner product $x\cdot y$ and norm $||x||=(x\cdot x)^{1/2}$, for vectors $x,y$ in $V$. Any finite sum $\sum\limits _{j=1}^k c_j \cdot x_j$, for real numbers $c_j$ and vectors $x_j$ in $V$, $1\leq j\leq k$, is called a linear combination of the vectors $\{ x_1 , x_2,\dots , x_k\}$.  A subset $E$ of $V$ is called linearly independent if for any finite subset $\{ x_1 , x_2 , \dots ,x_k\}$ of $E$ and real numbers $\{ c_1 , c_2 , \dots ,c_k\}$ the relation $\sum\limits _{j=1}^k c_j \cdot x_j=0$ implies that all $c_j=0$, $1\leq j\leq k$. We say that $E$ is pairwise orthogonal if $x\cdot y=0$, for any two different vectors $x$ and $y$ in $V$. As is well known, any pairwise orthogonal set is also linearly independent. The converse implication is in general false. A pairwise orthogonal set of vectors of unit length is called an ortonormal set, or an $ON-$set for short.\\[1em]
The linear span of a set $E\subset V$ is the set of all linear combinations $\sum\limits _{j=1}^k c_j \cdot x_j$ of vectors $x_j$ from $E$ and is denoted by $span(E)$. Clearly, $span(E)$ is a linear subspace of $V$.\\[1em]
We define the angle $\theta$ between two non--zero vectors $x$ and $y$ in $V$ by
\begin{displaymath} \cos \theta =\frac{x\cdot y}{||x||\cdot ||y||},\quad 0\leq \theta \leq \pi .\end{displaymath} 
If $x=0$ or $y=0$ we define the angle $\theta$ to be zero. In the following we identify tha angle $\theta$ with its $cosine$ $p$,  $p=\cos \theta$ and $-1\leq p\leq 1$. The Gram -- Schmidt Orthogonalization Algorithm produces an ortogoanal set of vectors from an linearly independent set. We state this result as follows, c.f. \cite{S}, Theorem 31, or \cite{KKOP}, \cite{RS}.
\begin{theorem} Let $E=\{ x_i\}$ be a finite or denumerably infinite linearly independent subset of $V$. Then there is a set $E^\prime =\{ y_i\}$ of vectors in $V$ with the following properties:
\begin{equation} span(\{ y_1 , y_2,\dots , y_i\} ) = span(\{ x_1 , x_2,\dots , x_i\} ),\, for \,  i=1, 2, \dots ,\end{equation}
\begin{equation} span(E^\prime)=span(E),\end{equation}
\begin{equation} \textrm{ the set } E^\prime \textrm{ is an $ON-$set}.\end{equation}\end{theorem}
We are going to generalize Theorem 1 in the following way. Fix any $-1<p<1$ and a linearly independent set $\{ x_1, x_2, \dots ,x_n\}$ in $V$. We construct a set of vectors $\{ z_1, z_2, \dots ,z_n\}$ that spans the same subspaces as $\{ x_1, x_2, \dots ,x_n\}$ and satisfies $||z_i||=1$, $z_i\cdot z_j=p$, for $1\leq i,j\leq n,$ $i\ne j$. It turns out that a necessary and sufficient condition for such a construction is that $p(n-1)+1>0$ in the $n-$dimensional case and $p\geq 0$ in the infinite dimensional case.
\begin{theorem} Let $\{ x_i\}_{i=1}^n$ be a  linearly independent subset of $V$ and let $-1<p<1$. Then there exists a set of vectors 
$\{ z_i\}_{n=1}^n$ such that
\begin{equation}span(\{ z_1 , z_2,\dots , z_i\} ) = span(\{ x_1 , x_2,\dots , x_i\} )\end{equation}
and 
\begin{equation}  ||z_i||=1, \quad z_i\cdot z_j=p, \end{equation}
for all $1\leq i,j\leq n$, $i\ne j$,  if and only if $-\frac{1}{n-1}<p<1$.
\end{theorem}
\begin{corollary} Let $\{ x_i\}_{n=1}^\infty$ be a denumerably infinite and linearly independent set of vectors in $V$ and let $-1<p<1$. Then there exists a set
$\{ z_i\}_{n=1}^\infty$ in $V$ such that (4) and (5) hold for all $1\leq i,j $, $i\ne j$, if and only if $0\leq p<1$.\end{corollary}
  Theorem 2 and Corollary 1 can be stated in terms of distances in stead of angles, since if $||x||=||y||=1$ then $||x-y||^2=2(1-x\cdot y)$.  
  \begin{theorem}
{\it Let $\{ x_i\}_{i=1}^n$ be a linearly independent set of vectors in $V$ and $0<d<2$. Then there exists vectors $\{ z_i\}_{i=1}^n$ satisfying (4) and 
\begin{equation}||z_i||=1\quad \textrm{and}\quad ||z_i-z_j||=d\end{equation}
for $1\leq i,j\leq n$, $i\ne j$, if and only if $0<d<\sqrt{\frac{2n}{n-1}}$.} \end{theorem}
\begin{corollary} Let $\{ x_i\}_{n=1}^\infty$ be a denumerably infinite and linearly independent set of vectors in $V$ and let $0<d<2$. Then there exists a set
$\{ z_i\}_{n=1}^\infty$ in $V$ such that (4) and (6) hold for all $1\leq i<j $, if and only if $0<d\leq \sqrt{2}$.\end{corollary}
{\it Remark 1.}  Taking $d=1$ in (6) we get an equilateral set $\{ 0,z_1,z_2,\dots ,z_n\}$ with $n+1$ points. For more on equilateral sets in normed spaces, see \cite{SV}.
\section{Proofs} This section contains the proofs of Theorem 2 and Corollary 1 and we begin with the proof of Theorem 2.\\[1em]
{\it Proof of Theorem 2.} Let $\{ x_i\} _{i=1}^n$ be as in the theorem, let $-\frac{1}{n-1}<p<1$ and perform the Gram--Schmidt algorithm on $\{ x_i\} _{i=1}^n$. This gives vectors $\{ y_i\} _{i=1}^n$ satisfying (1) and (3). We will define the vectors $\{ z_i\} _{i=1}^n$ inductively and start with $z_1=y_1$. Next let $w=y_2+a\cdot z_1$. Then $||w||^2=1+a^2$ and $w\cdot z_1=a$. We take $a$ such that $w\cdot z_1=p\cdot ||w||$, which gives
$a=p/\sqrt{1-p^2}$. Finally define $z_2=w/||w||$.\\[1em]
Now assume that we have defined $z_1, z_2, \dots ,z_k$ such that $||z_i||=1$ and $z_i\cdot z_j=p$, $1\leq i,j\leq k,$ $i\ne j$. Let $w=y_{k+1}-\sum \limits _{j=1}^k c_j\cdot z_j$. Then
\begin{displaymath} w\cdot z_i=-c_i  -p\cdot \sum\limits _{j\ne i}c_j,\quad 1\leq i\leq k,\end{displaymath}
and \begin{displaymath}||w||^2=1+||\sum\limits _{j=1}^k c_j\cdot z_j ||^2,\end{displaymath}
since $y_{k+1}$ is orthogonal to $span\{ z_1, z_2, \dots ,z_k\}\subset span \{ y_1, y_2, \dots ,y_k\}$. This gives
\begin{displaymath} w\cdot z_i=-c_i  -p\cdot \sum\limits _{j\ne i}c_j =p\cdot ||w||,\end{displaymath}
for $1\leq i\leq k$. Subtracting the last equations pairwise gives $c_1=c_2=\cdots  =c_k=d$ and 
\begin{equation}d  (1+p (k-1))=-p\cdot ||w||.\end{equation}
Further,
\begin{equation} ||w||^2=1+d^2\cdot ||\sum\limits _{i=1}^k z_i||^2=1+d^2(k+p (k^2-k) ).\end{equation}
Combining (7) and (8) gives
\begin{displaymath}d^2\cdot (1+p(k-1))\cdot (1-p)\cdot (1+pk)=p^2.\end{displaymath}
This equation can be solved for $d$, since $p>-\frac{1}{n-1}\geq -\frac{1}{k}.$ Now define $z_{k+1}=w/||w||$. Continuing in this way produces vectors $\{ z_i\} _{i=1}^n$ satisfying (5). It is easy to see that 
$\{ z_i\} _{i=1}^n$ is linearly independet, which proves (4).\\[1em]
Conversely, assume that $\{ z_1, z_2, \dots ,z_n\}$ are vectors in $V$ satisfying (4) and (5). Then $\{ z_1, z_2, \dots ,z_n\}$ is linearly independent and the formula
$0<||z_1+z_2+\cdots +z_n||^2=n(1+p(n-1))$ proves that $p>-1/(n-1)$.\hfill $\triangle$\\[1em]
{\it Proof of Corollary 1}. Let $\{ x_i\}_{i=1}^\infty $ be as in the corollary and $0\leq p<1$. The proof of Theorem 2 applies and gives vectors  
$\{ z_i\}_{i=1}^\infty $ with the desired properties.\\[1em]
Conversely, assume that $\{ z_i\}_{i=1}^\infty $ has properties (4) and (5) for $1\leq i<j.$ Then $\{ z_i\}_{i=1}^\infty $ is linearly independent and $0<||z_1+z_2+\cdots z_n||^2=n(1+p(n-1))$, $n\geq 1$. It follows that $p\geq 0$ by letting $n\rightarrow \infty$.\hfill $\triangle$
\section{The Gram--Schmidt p--Algorithm} In this section we describe the algorithm that produces the sequence $\{ z_i\}_{i=1}^n$ in Theorem 2.  Let $\{ y_i\}_{i=1}^n$ be an $ON-$set in $V$ and let  $-\frac{1}{n-1}<p<1$. We define the vectors $\{ z_i\}_{i=1}^n$ inductively by the method in the proof of Theorem 2. Put $z_1=y_1$ and assume that $z_2, z_2, \dots z_k$ have been defined. Then let
\begin{displaymath}w=\alpha \cdot y_{k+1}+\beta \cdot (z_1+z_2+\cdots +z_k),\end{displaymath}
where $\alpha$ and $\beta$ satisfy
\begin{displaymath} ||w||^2=\alpha ^2+\beta ^2(k+k(k-1)p)=1\end{displaymath}
and
\begin{displaymath} w\cdot z_i=\beta \cdot (1+(k-1)p)=p\cdot ||w||,\, 1\leq i\leq k.\end{displaymath}
An easy calculation yields
\begin{equation} \alpha ^2=\frac{(1-p)(1+pk)}{1+p(k-1)}\quad \textrm{and}\quad \beta =\frac{p}{1+p(k-1)},\end{equation}
which gives the following algorithm.\\[1em]
{\it The Gram--Schmidt p--Algorithm $(GS_p)$.} The input is a linearly independent set $\{ x_1 , x_2,\dots , x_n\}$ in a real inner product space $V$. The
output is a set $\{ z_1 , z_2,\dots , z_n\} $ of vectors with the properties (4) and (5) in Theorem 2.\\[1em]
{\it Step 0.} Perform the Gram--Schmidt Algorithm in Theorem 1 and get an $ON-$ set $\{ y_1 , y_2,\dots , y_n\}$.\\[1em]
{\it Step 1}. Define $z_1=y_1$.\\[1em]
{\it Step 2}. After vectors $\{ z_1 , z_2,\dots , z_k\} $ have been constructed, we define 
\begin{equation}z_{k+1}=\sqrt{\frac{(1-p)(1+pk)}{1+p(k-1)}}\cdot y_{k+1}+\frac{p}{1+p(k-1)}\cdot (z_1+z_2+\cdots +z_k),\end{equation}   \\[1em]
{\it Step 3}. The set $\{ z_1 , z_2,\dots , z_n\} $ that is constructed after {\it Step 2} has been performed $(n-1)$ times has the desired properties.\hfill $\triangle$
\\[1em]
  Formula (10) will allow us to get explicit expressions for each $z_k$. We start with the following lemma. \begin{lemma}Let $\{y_k\}_1^n$ be an ON--set in $V$ and let $\{ z_k\}_1^n$ be defined by ($GS_p$). Then there are numbers $z_k(i)$, $1\leq i\leq k$ and $1\leq k\leq n$, such that
 \begin{displaymath}z_k=z_k(1)y_1+z_k(2)y_2+\cdots z_k(k)y_k,\end{displaymath}
 where $z_{k+1}(i)=z_k(i)$, for $1\leq i\leq k-1$ and $1\leq k\leq n-1$.
\end{lemma}
{\it Proof.} The case $k=2$ follows from an explicit calculation. Assume that the lemma holds for some fixed $2\leq k\leq n-1$ and all $1\leq i\leq k-1$ Then for $1\leq i\leq k$ by (10)
\begin{displaymath}z_{k+1}(i)=\frac{p}{1+p(k-1)}\cdot \left( z_1(i)+\cdots + z_k(i)\right)\end{displaymath}
and
\begin{displaymath}z_{k+2}(i)=\frac{1}{1+pk}\cdot (z_1(i)+\cdots +z_k(i)+z_{k+1}(i)),\end{displaymath}
from which the lemma follows.\hfill $\triangle$\\[1em]
We can now give explicit formulas for $z_k$, $1\leq k\leq n$. Define
\begin{equation}a_m(p)=\frac{p}{1+p(m-1)}\cdot \sqrt{\frac{(1-p)(1+p(m-1))}{1+p(m-2)}},
\end{equation}
m=1,2,\dots . Then we get the following result.
\begin{theorem} Let $\{y_k\}_1^n$ be an ON--set in $V$ and let $\{ z_k\}_1^n$ be defined by $(GS_p)$. Then $z_1=y_1$ and $z_k$ is given by
\begin{equation}z_k=a_1(p)y_1+a_2(p)y_2+\dots +a_{k-1}(p)y_{k-1}+\frac {1+p(k-1)}{p}\cdot a_k(p)y_k,\end{equation}
for $2\leq k\leq n$.\end{theorem}
{\it Proof}. The case $k=2$ follows by calculation. Assume that the theorem holds for some $2\leq k\leq n-1$ and consider $z_{k+1}$. We note that $z_{k+1}(i)$, $1\leq i\leq k-1$, and $z_{k+1}(k+1)$ are correct by Lemma 1 and (10). Finally, 
\begin{displaymath}z_{k+1}(k)=\frac{p}{1+p(k-1)}\cdot z_k(k)=a_k(p),\end{displaymath}
by (10) and (11), which completes the proof of the theorem\hfill $\triangle$
\section{Asymptotics.}
The final section of this paper is about the infinite dimensional case. Let $V$ be the standard real vector space $l^2$ of points $x=(x_1,x_2,\dots ,x_n,\dots )$ with finite norm $||x||=(x_1^2+x_2^2+\cdots +x_n^2+\cdots )^{1/2}$ and inner product $x\cdot y=x_1y_1+x_2y_2+\cdots +x_ny_n+\cdots .$ Let $\{ y_i\}_1^\infty$ be a denumerably infinite ON--set in $V$, $0<p<1$, and let $\{ z_k\}_1^\infty$ be the vectors obtained from Theorem 4 by identifying the finite dimensional space $V$ as a subspace of $l^2$ in a natural way. Then
\begin{equation} ||z_k||=1,\quad z_k\cdot z_m=p,\quad \textrm{and}\quad ||z_k-z_m||= \sqrt{2-2p},\end{equation} for $1\leq k<m$, in analogy with Theorems 2 and 3. Define a vector $z_0$ by 
\begin{displaymath}z_0=a_1(p)y_1+a_2(p)y_2+\dots +a_m(p)y_m+\dots .\end{displaymath} 
It is easy to see that $z_0$ is in $V$ and in fact that $|| z_0||=\sqrt{p}$ by a telescoping argument. We can now use $z_0$ to get the following asymptotic formula for $z_k$, as $k\rightarrow\infty$.
\begin{theorem} Let $\{ z_k\} _1^\infty$ be as above, then
$z_k\approx z_0+\sqrt{1-p}\cdot y_k,$ as $k\rightarrow \infty$, in the sense that
\begin{displaymath}\limsup\limits _{k\rightarrow\infty}\sqrt{k}\cdot  ||z_k-z_0-\sqrt{1-p}\cdot y_k||\leq \sqrt{p(1-p)}.
\end{displaymath}
\end{theorem}
{\it Proof.} Clearly, $||z_k-z_0-\sqrt{1-p}\cdot y_k||$ is at most
\begin{displaymath}  \left( \sum\limits _k^\infty |a_m(p)|^2\right) ^{1/2}+
|\sqrt{\frac{(1-p)(1+p(k-1))}{1+p(k-2)} }-\sqrt{1-p}|=A+B.\end{displaymath}
We get
\begin{displaymath} A^2=\sum\limits _k^\infty |a_m(p)|^2=\frac{p^2(1-p)}{1+p(k-2)}\end{displaymath}
by a telescoping series and 
\begin{displaymath}B\leq \frac{1}{2}\cdot \frac{p\sqrt{1-p}}{1+p(k-2)}\end{displaymath}
 by an elementary calculation. The theorem follows easily.\hfill $\triangle$
  \\[1em] {\it Remark 2.}  The set $\{ z_i\} _1^\infty$ in (13) is an infinite equilateral set in $V$ for $p=1/2$.\\[1em] {\it Historical remark.} J\o rgen Pedersen Gram, 1850 -- 1916, was a danish mathematician working in booth pure and applied mathematics. The german mathematician Erhard Schmidt, 1876 -- 1959, is considered as one of the fathers of modern abstract functional analysis. The Gram--Schmidt algorithm was however presented before both Gram and Schmidt by the famous french mathematician Pierre -- Simon Laplace, 1749 -- 1827. 

Address: Department of Mathematics and Mathematical Statistics, University of Ume\aa , S--901 87 
Ume\aa , Sweden. E--mail: tord.sjodin@math.umu.se
 
\end{document}